\documentclass{article}
\usepackage{latexsym}
\usepackage{amssymb}
\usepackage{amsfonts}
\usepackage{amsmath}
\usepackage{graphicx}
\usepackage{amsfonts}
\usepackage{epstopdf}
\usepackage{color}
\usepackage{enumerate}
\usepackage[latin1]{inputenc}
\usepackage{mathtools}
\usepackage{amsthm}
\usepackage[latin1]{inputenc}
\usepackage{mathrsfs}
\usepackage{graphicx}
\usepackage[latin1]{inputenc}

\newcommand{\FF}{{\cal F}}

\def\lf{\left\lfloor}
\def\rf{\right\rfloor}

\newtheorem{teo}{Theorem}[section]

\begin{document}
\title{Rate of convergence of uniform transport processes  to Brownian sheet}
\date{}
\author{Carles Rovira\footnote{C. Rovira is supported by the grant PGC2018-097848-B-I00.}}

\maketitle

$\mbox{ }$\hspace{0.1cm} {\rm Departament de Matem\`atiques i Inform\`atica,
Universitat de Barcelona, Gran Via 585, 08007 Barcelona}. {\tt carles.rovira@ub.edu}

\begin{abstract}
In a previous paper we have constructed a family of processes, starting from a set of independent standard Poisson processes, that has realizations that converge almost surely to the Brownian sheet, uniformly in the unit square. Now, a rate of convergence from these processes to Brownian sheet is given.

\end{abstract}

\section{Introduction}

Let $W=\{W(s,t): (s,t)\in [0,1]^{2}\}$ be a Brownian sheet, i.e. a zero mean real
continuous Gaussian process with covariance function
$E[W{(s_1,t_1)} W{(s_2,t_2)}]=(s_1\wedge s_2)(t_1\wedge t_2)$ for any
$(s_1,t_1), (s_2,t_2)\in [0,1]^2$.
Our aim is to obtain the rate of convergence of strong approximations of the Brownian sheet.
This result is not only interesting from a purely mathematical point of view, but are
of great interest in order to provide sound approximation strategies to solutions of
stochastic partial differential equations whih arise in many fields as physics, biology or finance.

The study of the approximations of the Brownian sheet by uniform transport processes or processes constructed from a Poisson process  begins with the proof of the
weak convergence.
Bardina and Jolis \cite{BJ} prove that the process
\begin{equation*}
\frac{1}{n} \int_0^{tn}  \int_0^{sn} \sqrt{x y} (-1)^{N(x,y)}dxdy,
\end{equation*}
where $\{N(x,y), x\ge0, y \ge 0 \}$ is a Poisson process in the plane, converges in
law to a Brownian sheet when $n$ goes to infinity
and Bardina, Jolis and Rovira \cite{BJR} extended this result to  the $d$--parameter Wiener processes.

On the other hand,  Bardina, Ferrante and Rovira \cite{BFR} constructed a family of processes, starting from a set of independent standard Poisson processes, that has realizations that convergence almost surely to a Brownian sheet.
Our purpose is to give the rate of convergence of such approximations.
As far as we know, our work is the first rate of convergence for the multiparameter case for
this family of approximations..


There exist several literature about strong convergence of uniform transport processes and the study of the corresponding rate of convergence. In the seminal paper of Griego, Heath and Ruiz-Moncayo~\cite{art G-H-RM}, the authors
presented
realizations of a sequence of the uniform transform processes that converges almost surely
to the standard Brownian motion, uniformly on the unit time interval.
In \cite{art G-G2} Gorostiza and Griego
extended the result of \cite{BBR} to the case of diffusions. Again Gorostiza and Griego
\cite{art G-G} and Cs\"{o}rg\H{o} and Horv\'ath \cite{CH} obtained the
rate of convergence of the approximation sequence.
More recently, Garz\'on, Gorostiza and Le\'on \cite{GGL} defined a sequence of
processes that converges strongly to fractional Brownian motion
uniformly on bounded intervals, for any Hurst parameter $H\in(0,1)$
and computed the rate of convergence. In \cite{GGL2} and \cite{GGL3}
the same authors deal with subfractional Brownian motion and
fractional stochastic differential equations.
Bardina, Binotto and Rovira \cite{BBR} proved the strong convergence to
a complex Brownian motion and obtained the corresponding rate of convergence.

The structure of the paper is the following. In the next section we recall the approximations that
converge almost surely to the Brownian sheet and we present our theorem. In the last section we give the proof of our result. It is based on a combination of the properties of the Brownian sheet and the use of the rate of convergence for the Brownian motion given by 
Griego, Heath and Ruiz-Moncayo in \cite{art G-H-RM}.

\section{Approximations and main result}

Let us recall the approximation processes introduced in \cite{BFR}.
For any $n$ and $\lambda >0$,  consider the partition of  the unit square
$[0,1]^2$ in disjoint rectangles
\[
\Big([0,\frac{1}{n^\lambda}] \times[0,1]\Big) \cup
\Big(  \bigcup_{k=2}^{\lf n^\lambda\rf}  (\frac{k-1}{n^\lambda} , \frac{k}{n^\lambda}] \times[0,1]  \Big)
\cup
\Big((\frac{\lf n^\lambda\rf}{n^\lambda} , 1] \times[0,1]\Big) .
\]
where $\lf x \rf$ denotes the greatest integer less than or equal to $x$.

If $W=\{W(s,t): $ $\,(s,t)\in [0,1]^{2}\}$ is a Brownian sheet on the unit square,
let $W^k$ denotes its restriction to each of the above defined rectangles
$ (\frac{k-1}{n^\lambda} , \frac{k}{n^\lambda}] \times [0,1]$.
That is, 
\[
W^k(t):=W(\frac k{n^\lambda},t)-W(\frac{k-1}{n^\lambda},t),
\]
for $k\in\{1,2,\dots,\lf n^\lambda\,\rf\}$.
Thus, for any $l \in\{1,2,\dots,\lf n^\lambda\,\rf\}$ and $t \in [0,1]$
\[
W(\frac{l}{n^\lambda},t) = \sum_{k=1}^l W^k(t).
\]
Moreover, putting $\tilde W^k(t):=n^{\frac{\lambda}{2}} W^k(t)$, we obtain a
family
\[
\{\tilde W^k;\,k\in\{1,2,\ldots \lf n^\lambda\,\rf\} \}
\]
of
independent standard  Brownian motions defined in $[0,1]$.

From the paper of Griego, Heath and Ruiz-Moncayo  \cite{art G-H-RM} it is  known that
there exist realizations of  uniform transport processes
that converge strongly and uniformly on bounded time intervals
to Brownian motion.
So, we can get an approximation sequence $\{ \tilde W^{(n)k};\,n \ge 1\} $ for each one of the standard Brownian motions
$\tilde W^k;\,k\in\{1,2,\dots,\lf n^\lambda\,\rf\}$. We can state Theorem 1 in \cite{art G-G} for such approximation sequence for any $k$:
\begin{teo}\label{ratebrownian} 
	There exists a version  $\{ \tilde W^{(n)k}(t), t \ge 0 \}$ of the uniform transport processes on the  same probability space as a Brownian motion process  $\{ \tilde W^k(t), t \ge 0 \}$, $\tilde W^k(0)=0$ so that
	$$\lim_{n \to \infty} \max_{0 \le t \le 1} |\tilde W^{(n)k}(t)
-\tilde W^k(t) |=, \qquad a.s.,$$
and such that for all $q > 0$
$$
P \big( \max_{0 \le t \le 1} |\tilde W^{(n)k}(t)
-\tilde W^k(t) | > \alpha n^{- \frac12} (\log n )^{\frac52} \big) = o(n^{-q}), \quad as \quad n \to \infty,$$where $\alpha$ is a positive constant depending on $q$.
\end{teo}

Then,  the Brownian sheet is approximated by a process
$W_n$ such that for any   $l \in\{1,2,\dots,[n^\lambda\,]\}$ and $t \in [0,1]$
$$W_n(\frac{l}{n^\lambda},t) =\sum_{k=1}^l W^{(n)k}(t)
= \sum_{k=1}^l \frac{1}{n^\frac{1+\lambda} 2} (-1)^{A_k}\int_0^{nt}(-1)^{N_k(u)}du,
$$
where
$$ W^{(n)k}(t)=\frac1{n^\frac{\lambda}{2}}\tilde W^{(n)k}(t),$$
and $\{N_k, k \ge 1\}$ is a family of independent standard Poisson processes and  $\{A_k, k \ge 1\}$ is a
sequence of independent random variables with law
$\textrm{Bernoulli}\left(\frac12\right)$, independent of the Poisson
processes.
Using linear interpolation, define $W_n(s,t)$
on the whole unit square as follows:
\begin{eqnarray}
W_{n} (s,t)&=&W_n (\frac{\lf s n^\lambda\rf}{n^\lambda} ,t) +(s n^\lambda -\lf s n^\lambda\rf)W_n (\frac{\lf s n^\lambda\rf+1}{n^\lambda} ,t),
\nonumber
\end{eqnarray}
for any $(s,t) \in [0,\frac{\lf n^\lambda\rf}{n^\lambda} ]\times [0,1]$ and
$W_{n} (s,t)=W_{n} (\frac{\lf n^\lambda\rf}{n^\lambda},t)$ for any
$(s,t) \in [\frac{\lf n^\lambda\rf}{n^\lambda},1 ]\times [0,1]$.

\bigskip

In the following theorem we give our main result, the rate of convergence of these processes:

\begin{teo}\label{ratesheet} 
		There exists realizations of the process $\{ W_{n} (s,t), (s,t) \in [0,1]^2 \}$ with $\lambda \in (0,\frac15)$
	on the same probability space as a  Brownian sheet $\{ W (s,t), (s,t) \in [0,1]^2 \}$
	such that
	$$ \lim_{n\rightarrow \infty}
	\max_{0\leq s,t\leq1}
	|W_n(s,t)-W(s,t)|=0 \quad a.s. $$
	and such that for all $\beta < \frac{\lambda}{2}$ and  $q > 0$
	$$
	P \big( 	\max_{0\leq s,t\leq1}
	|W_n(s,t)-W(s,t)|> \alpha n^{- \beta}  \big) = o(n^{-q}), \quad as \quad n \to \infty,$$where $\alpha$ is a positive constant depending on $q$.
\end{teo}

Notice that the rate of convergence if worse than the rate for the Brownian motion case due to the properties of the Brownian sheet.

The first part of the Theorem has been proved in Theorem 2.1 in \cite{BFR}. In the next section we give the proof of the rate of convergence.

\section{Proof}

We will begin recalling some technical results about submartingales.
Let us begin with a version of Theorem 1, page 74  from Imkeller \cite{I}. 

\begin{teo}\label{imkeller}
Let $M$ be a non-negative submartingale and set $\psi (t)= t \log^+ t,$ for any $t\ge 0$. Then, there is a constant $c$ which does not depend on $M$ such that, for any $\beta>0$
$$
P \Big( \max_{(s,t) \in [0,1]^2} | M(s,t) | > \beta \Big) \le\frac{c}{\beta} \Vert M(1,1) \Vert_\psi,$$	
where
$$\Vert X\Vert_\psi= \inf \{ \mu >0: E \Big( \psi \big( \frac{|X|}{\mu}\big)\Big) \le 1  \}.$$
	\end{teo}

We will apply this theorem to the process $$M(s,t):=\exp(B(s,t)),\quad  (s,t) \in [0,1]^2$$ where $\{B(s,t): (s,t)\in [0,1]^{2}\}$ is a Brownian sheet. Considering the family of $\sigma$-algebras $\FF_{s,t}=\sigma <B(u,v); (u,v) \in [0,s]\times[0,t] >$, $M(s,t)$ is clearly nonnegative, integrable, $\FF_{s,t}$-measurable and it is a submartingale. Indeed, if $s'<s$ and $t'<t$
\begin{eqnarray*}
& &E(e^{B(s,t)}|\FF_{s',t'})=
E(e^{B(s,t)-B(s',t')}e^{B(s',t')}|\FF_{s',t'}) \\
& &\quad = E(e^{B(s,t)-B(s',t')})e^{B(s',t')}
=e^{\frac{tt'-ss'}{2}} e^{B(s',t')} \ge e^{B(s',t')}.
\end{eqnarray*}

On the other hand,  $\Vert  M(1,1) \Vert_\psi = \Vert \exp(B(1,1)) \Vert_\psi$ is a finite constant. 
Indeed, using that $B(1,1)$ is a $N(0,1)$, 
\begin{eqnarray*}
& &E\Big( \frac{\exp(B(1,1))}{\mu} 
\log^+ \big(\frac{\exp(B(1,1))}{\mu} 
\big)\Big)= \frac{1}{(2\pi)^\frac12}\int_{\log \mu}^\infty \frac{e^x}{\mu} (x - \log \mu)  e^{-\frac12 x^2} dx\\
& &\qquad=\frac{1}{(2\pi)^\frac12} \int_{0}^\infty \frac{e^{y+\log \mu}}{\mu} y e^{-\frac12 (y+\log \mu)^2} dx\\
& &\qquad= \frac{1}{(2\pi)^\frac12} \frac{e^{\frac12} }{\mu} \int_{0}^\infty  y e^{-\frac12 (y+(\log \mu-1))^2} dx\\
& &\qquad= \frac{1}{(2\pi)^\frac12} \frac{e^{\frac12 }}{\mu} \big(  e^{-\frac12 (\log \mu-1)^2} + \big( \frac{\pi}{2} \big)^\frac12 (1- \log \mu)\big),
\end{eqnarray*}
and this quantity is bigger that 1 for $\mu$ sufficiently small and smaller that 1 for $\mu$ big enough. We have used that  $$ \int_{0}^\infty  x e^{-\frac12 (x-m)^2} dx =   e^{-\frac12 m^2} + m \big( \frac{\pi}{2} \big)^\frac12. $$

\medskip

Let us give now the proof of our Theorem:

\medskip

{\it Proof of Theorem \ref{ratesheet}:}
Our aim is to bound
$$
P:=P \Big( \max_{(s,t) \in [0,1]^2}
|W_n(s,t)-W(s,t)| >K_n \Big )
$$
where $K_n=c n^{-\beta}$, with $\beta \in (0,\frac{\lambda}{2} )$.

Fixed $s \in [0,1]$, $n \ge 1$, there exists $k \in \{0,\ldots,\lf n^\lambda\rf\}$ such that $$\frac{k}{n^\lambda} \le s < \frac{k+1}{n^\lambda},$$
that is $k=\lf s n^\lambda\rf$. 
So, we can write

\begin{eqnarray*}
P  	&\le &P \Big(\max_{(s,t) \in [0,1]^2}
	|W_n(s,t)-W_n(\frac{\lf s n^\lambda\rf}{n^\lambda},t)| > \frac{K_n}{3} \Big )
	\\ &+ &
P \Big(\max_{(s,t) \in [0,1]^2}
|W_n(\frac{\lf s n^\lambda\rf}{n^\lambda},t)-W(\frac{\lf s n^\lambda\rf}{n^\lambda},t)| > \frac{K_n}{3} \Big )
\\&+ &
P \Big(\max_{(s,t) \in [0,1]^2}
|W(\frac{\lf s n^\lambda\rf}{n^\lambda},t)- W(s,t)| > \frac{K_n}{3} \Big )
\\&:=& P_1 + P_2+P_3.
\end{eqnarray*}

Let us study first $P_2$. Using the definitions of $W^{(n)k}$ and $ W^k$, for any $k$, 
we get
\begin{eqnarray*}
P_2	&=&P \Big( \max_{1\leq l\leq  \lf n^\lambda\rf} \max_{ t \in [0,1]}
	|\sum_{k=1}^l W^{(n)k}(t)-\sum_{k=1}^l W^k(t)| > \frac{K_n}{3}  \Big )
	\\
		&=&P \Big( \max_{1\leq l\leq  \lf n^\lambda\rf} \max_{ t \in [0,1]}
	|\sum_{k=1}^l \tilde W^{(n)k}(t)-\sum_{k=1}^l \tilde W^k(t)| > n^{\frac{\lambda}{2}}\frac{K_n}{3}  \Big )
	\\
	&\le &P \Big( 
	\sum_{k=1}^{\lf n^\lambda\rf}  \max_{ t \in [0,1]}| \tilde W^{(n)k}(t)- \tilde W^k(t)| > n^{\frac{\lambda}{2}}\frac{K_n}{3}  \Big )
	\\
		&\le & 	\sum_{k=1}^{\lf n^\lambda\rf}  P \Big( 
 \max_{ t \in [0,1]}| \tilde W^{(n)k}(t)- \tilde W^k(t)| > n^{-\frac{\lambda}{2}}\frac{K_n}{3}  \Big )
.
\end{eqnarray*}
Finally, from the rate of convergence from $W^{(n)k}$ and $\tilde W^k$, for any $k$,  Theorem \ref{ratebrownian},
we get
\begin{equation*}
	P_2 \le 
	n^\lambda  P \Big( 
	\max_{ t \in [0,1]}| \tilde W^{(n)1}(t)- \tilde W^1(t)| > n^{-\frac{\lambda}{2}}\frac{K_n}{3}  \Big ) =o(n^{-q}),
\end{equation*}
for any $q>0$. 

Let us consider now $P_3$. We can write, using the well-known properties of the Brownian sheet 
\begin{eqnarray*}
	P_3	&=&P \Big( \max_{0\leq l\leq  \lf n^\lambda\rf} \max_{ t \in [0,1]}
	\max_{ r \in [0,n^{-\lambda}]}
	| W(\frac{l}{n^\lambda}+r,t)-W(\frac{l}{n^\lambda},t)| > \frac{K_n}{3}  \Big )
	\\
	&\le& \sum_{l=0}^{\lf n^\lambda\rf}  P \Big(  \max_{ t \in [0,1]}
	\max_{ r \in [0,n^{-\lambda}]}
	| W(\frac{l}{n^\lambda}+r,t)-W(\frac{l}{n^\lambda},t)| > \frac{K_n}{3}  \Big )
	\\
		&\le& (n^\lambda+1)  P \Big(  \max_{ t \in [0,1]}
	\max_{ r \in [0,n^{-\lambda}]}
	| W(r,t)| > \frac{K_n}{3}  \Big )
	\\
	&=&  (n^\lambda+1)P \Big(  \max_{ (s,t) \in [0,1]^2}
| W(\frac{s}{n^\lambda},t)| > \frac{K_n}{3}  \Big ).
\end{eqnarray*}
We have considered that $W(\frac{l}{n^\lambda}+r,t)=W(1,t)$ if $\frac{l}{n^\lambda}+r>1$.
Let us define $B(s,t):= n^\frac{\lambda}{2} W(\frac{s}{n^\lambda},t)$, for any $(s,t) \in [0,1]^2$. Clearly $B$ is a Brownian sheet and we can bound $P_3$ as
$$
P_3 \le  (n^\lambda+1) P \Big(  \max_{ (s,t) \in [0,1]^2}
| B(s,t)| > \frac{K_n}{3} n^\frac{\lambda}{2} \Big ).$$
So, using Theorem \ref{imkeller}
\begin{eqnarray*}
P_3 & \le & 2 (n^\lambda+1) P \Big(  \max_{ (s,t) \in [0,1]^2}
 B(s,t) > \frac{K_n}{3} n^\frac{\lambda}{2} \Big )\\
 & =& 2 (n^\lambda+1) P \Big(  \max_{ (s,t) \in [0,1]^2}
 \exp(B(s,t)) > \exp(\frac{K_n}{3} n^\frac{\lambda}{2} )\Big ) 
 \\
 & =& 2 (n^\lambda+1) c \exp(- \frac{K_n}{3} n^\frac{\lambda}{2} ) \Vert \exp(B(1,1)) \Vert_\psi \\
 & =& 2 (n^\lambda+1) c_2 \exp(- \frac{K_n}{3} n^\frac{\lambda}{2} )
 \\
 & =& 2 (n^\lambda+1) c_2 \exp(- c_3 n^{\frac{\lambda}{2} - \beta} ),
 \end{eqnarray*}
since $\Vert \exp(B(1,1)) \Vert_\psi$ is a finite positive constant that does not depend on $n$. 

Finally, let us study $P_1$. Using that we have defined the intermediate points using linear interpolation we have that for any $s \in [\frac{l}{n^\lambda},\frac{l+1}{n^\lambda} )$ we get that
$$
|W_n(s,t)-W_n(\frac{\lf s n^\lambda\rf}{n^\lambda},t)| \le 
|W_n(\frac{l}{n^\lambda},t)-W_n(\frac{l+1}{n^\lambda},t)|$$
and we can write that

\begin{eqnarray*}
	P_1	&=&P \Big( \max_{0\leq l\leq  \lf n^\lambda\rf} \max_{ t \in [0,1]}
	|W_n(\frac{l}{n^\lambda},t)-W_n(\frac{l+1}{n^\lambda},t)| > \frac{K_n}{3}  \Big )
	\\
	&\le&P \Big( \max_{0\leq l\leq  \lf n^\lambda\rf} \max_{ t \in [0,1]}
	|W_n(\frac{l}{n^\lambda},t)-W(\frac{l}{n^\lambda},t)| > \frac{K_n}{9}  \Big )
	\\
	&&+ P \Big( \max_{0\leq l\leq  \lf n^\lambda\rf} \max_{ t \in [0,1]}
	|W(\frac{l}{n^\lambda},t)-W(\frac{l+1}{n^\lambda},t)| > \frac{K_n}{9}  \Big )
	\\
	&&+ P \Big( \max_{0\leq l\leq  \lf n^\lambda\rf-1} \max_{ t \in [0,1]}
	|W(\frac{l+1}{n^\lambda},t)-W_n(\frac{l+1}{n^\lambda},t)| > \frac{K_n}{9}  \Big )\\
	&:=& P_{1,1}+P_{1,2}+P_{1,3},
\end{eqnarray*}
where if $\frac{l+1}{n^\lambda} > 1$ we assume that $W(\frac{l+1}{n^\lambda},t)=W_n(\frac{l+1}{n^\lambda},t)=W(1,t)$.

Notice that $P_{1,2}$ can be bounded as $P_3$ and on the other hand, $P_{1,1}$ and $P_{1,3}$ can be studied as $P_2$. Putting together the bounds for $P_1$, $P_2$ and $P_3$ we finish the proof.

\end{document}